\documentclass[12pt,a4paper,fleqn]{article}
\usepackage{amsmath}
\usepackage{amssymb}
\usepackage{amscd}
\usepackage[dvips]{graphicx}

\newtheorem{theorem}{Theorem}

\newtheorem{corollary}[theorem]{Corollary}
\newtheorem{definition}[theorem]{Definition}
\newtheorem{example}[theorem]{Example}
\newtheorem{lemma}[theorem]{Lemma}
\newtheorem{proposition}[theorem]{Proposition}
\newtheorem{remark}[theorem]{Remark}
\newtheorem{conjecture}[theorem]{Conjecture}

\newenvironment{proof}{{\bf Proof. }}{\hfill$\rule{1ex}{1ex}$\par\medskip}

\begin{document}
\newcommand{\bt}{\begin{theorem}}
\newcommand{\et}{\end{theorem}}
\newcommand{\bd}{\begin{definition}}
\newcommand{\ed}{\end{definition}}
\newcommand{\bs}{\begin{proposition}}
\newcommand{\es}{\end{proposition}}
\newcommand{\bp}{\begin{proof}}
\newcommand{\ep}{\end{proof}}
\newcommand{\be}{\begin{equation}}
\newcommand{\ee}{\end{equation}}
\newcommand{\ul}{\underline}
\newcommand{\br}{\begin{remark}}
\newcommand{\er}{\end{remark}}
\newcommand{\bex}{\begin{example}}
\newcommand{\eex}{\end{example}}
\newcommand{\bc}{\begin{corollary}}
\newcommand{\ec}{\end{corollary}}
\newcommand{\bl}{\begin{lemma}}
\newcommand{\el}{\end{lemma}}
\newcommand{\bj}{\begin{conjecture}}
\newcommand{\ej}{\end{conjecture}}

\def\beal{\begin{array}{l}}
\def\beac{\begin{array}{c}}
\def\beacl{\begin{array}{cl}}
\def\ena{\end{array}}

\def\diy{\displaystyle}

\title{Nilpotent Steiner loops of class $2$}

\author{A. Grishkov, D. Rasskazova, M. Rasskazova, I. Stuhl}

\date{}
\footnotetext{\par}

\maketitle

\begin{abstract}
We describe Steiner loops of nilpotency class $2$ and establish the classification of finite 3-generated nilpotent Steiner loops of nilpotency class $2$.
\end{abstract}

\section{Introduction and Preliminaries}

A {\it loop} is a set $L$ with a binary operation $\cdot$ and a neutral element $1 \in L$, such that for every $a, b \in L$ the equations $a\cdot x = b$ and $y\cdot a = b$ have unique solutions in $L$.

The {\it center} of a loop $L$ is an associative subloop $Z(L)$ consisting of all elements of
$L$ which commute and
associate with all other elements of $L$.
A loop $L$ is {\it nilpotent} if the series $L$, $L/Z(L)$,
$[L/Z(L)]/Z[L/Z(L)]$ ... terminates at $1$ in finitely many steps.
In particular, $L$ is of {\it nilpotency class two} if $L/Z(L)\neq 1$ and is an abelian group.
For $x, y, z \in L$, define the {\it associator} $(x, y, z)$ of $x, y,
z$ by $(xy)z = (x(yz))(x, y, z)$. The {\it associator subloop}
$Ass(L)$ of $L$ is the smallest normal subloop $H$ of $L$ such that
$L/H$ is a group. Thus, $Ass(L)$ is the smallest normal subloop
of $L$ containing associators $(x, y, z)$ for all $x, y, z\in L$.

A {\it Steiner triple system} is an incidence
structure consisting of points and blocks such that every two
distinct points are contained in precisely one block and any block
has precisely three points. A finite Steiner triple system with
$n$ points exists if and only if $n \equiv 1$ or $3$ $(mod\hskip
4pt 6)$ (Cf. \cite{CR}).

A Steiner triple system $\mathfrak{S}$
induces a multiplication on pairs of distinct points $x, y$ taking as the product the
third point of the block joining $x$ and $y$.
Defining $x\cdot x=x$, we get a {\it Steiner quasigroup}
associated with $\mathfrak{S}$. Adjoining an element $e$ with $ex=xe=x$, $xx=e$, we obtain
a Steiner loop $S$, see below.
Conversely, a Steiner loop determines a Steiner triple
system whose points are the elements of $S\setminus\{e\}$, and
the blocks are triples $\{x, y, xy\}$ where $x, y\in \mathfrak{S}$, $x\neq y$.

A loop $L$ is said to be {\it totally symmetric} if $x\cdot y=y\cdot x$ and $x\cdot(x\cdot y)=y$ for all $x, y\in L$.
A totally symmetric loop of exponent $2$ is called a {\it Steiner loop}. A variety of universal algebras is called a {\it Schreier variety} if every subalgebra of any free algebra in that variety is also free in that variety.
Steiner loops form a Schreier variety; it is precisely the variety of all
diassociative loops of exponent $2$ (see in \cite{GP} p. 310). Steiner loops are in a one-to-one correspondence with Steiner triple systems. Multiplication groups of Steiner loops have been studied in \cite{SS}. Central extensions of Steiner loops and quasigroups yielding algebraic structures of Steiner triple systems with cyclic orientations on each triple have been introduced in \cite{SS2, S}.

Since Steiner loops form a variety, we can deal with free objects and according to
\cite{GRRS}, we have the following.
Let $\tt X$ be a finite ordered set and let $W(\tt X)$ be a set of non-associative $\tt X$-words.
The set $W(\tt X)$ has an order, $>$, such that
$v>w$ if and only if $|v|>|w|$ or $|v|=|w|>1$, $v=v_1v_2,$ $w=w_1w_2,$ $v_1>w_1$ or
$v_1=w_1,v_2>w_2$. Next, we define
the set $S({\tt X})^{\ast}\subset W({\tt X})$ of $S$-words by induction upon the length of the word:
\begin{itemize}
  \item ${\tt X}\subset S({\tt X})^{\ast}$,
  \item $wv\in S(\tt X)^{\ast}$ precisely if, $v, w\in S(\tt X)^{\ast}$, $w> v$ and if
  $w=w_1w_2$, then $v\neq w_i$, $(i=1, 2)$.
\end{itemize}
\noindent
On $S({\tt X})=S({\tt X})^{\ast} \cup \{\emptyset\}$ we define a multiplication in the following manner:
\begin{enumerate}
  \item $v\cdot w=w\cdot v=vw$ if $vw\in S({\tt X})$,
  \item $(vw)\cdot w=w\cdot (vw)=w\cdot(wv)=(wv)\cdot w=v$,
  \item $v\cdot v=\emptyset$.
\end{enumerate}
A word $v(x_1, x_2, ... ,x_n)$ is {\it irreducible} if $v\in S({\tt X})^{\ast}$.
The set $S({\tt X})$ with the multiplication as above is a free Steiner loop with the set of free generators $\tt X$.

In what follows we discuss Steiner loops of nilpotency class $2$.

\section{Centrally nilpotent Steiner loops of class 2}

Let $S({\tt X})>S_1({\tt X})>S_2({\tt X})>...$ be a central series of the free Steiner loop $S({\tt X})$ with free generators
${\tt X}=\{x_1,...,x_n\}$. Then $V=S({\tt X})/S_1({\tt X})$ is an $\mathbf{F}_2$-space
of dimension $n:=|{\tt X}|$. Given $\sigma =\{i_1>i_2>...>i_s\}\subseteq I_n$ define the corresponding element
$\sigma=(((x_{i_2}x_{i_1})x_{i_3})... x_{i_s})$ of $S({\tt X})$. As $S({\tt X})/S_1({\tt X})$ is an abelian $2$-group, it is isomorphic to $\mathbf{F}_2^n$. Hence, to any element $v\in \mathbf{F}_{2}^{n}$ an element $\sigma=(i,j,...)$ can be related, where $i, j, ...$ are the numbers of coordinates having value $1$ of the vector $v$. Therefore, $\{\sigma | \sigma \subseteq I_n\}$
is a set of representatives of $S({\tt X})/S_1({\tt X})$. Determine a set of representatives of
$Z=S_1({\tt X})/S_2({\tt X})$.

Set $L_f$ to be a central extension of $\mathbf{F}_2$-spaces $V$ and $Z$ in the variety of Steiner loops.
It is well known that $L_f$ is a central extension of
$Z$ by $V$ if and only if $L_f$ is isomorphic to a loop defined on
$V \times Z$ by the multiplication
\begin{equation}\label{muvelet}
(v_1, z_1)\circ (v_2, z_2) = (v_1 + v_2, f(v_1, v_2)+z_1+z_2).
\end{equation}
Here $f : V \times V \longrightarrow Z$ is a {\it Steiner loop cocycle}, that is, a map satisfying
\begin{equation}\label{cocycle}
f(0_V,v_1)=f(v_1,v_1)=0_Z, \hskip 2pt f(v_1,v_2)=f(v_2,v_1), \hskip 2pt
f(v_1 + v_2, v_2)=f(v_1,v_2)
\end{equation}
for all $v_1, v_2\in V$.
Denote by $\mathcal{Z}^2(V,Z)$ the set of all Steiner loop cocycles. Next, let $\mathcal{C}^1(V,Z)$ be the set of all functions $g:V \longrightarrow Z$ and
$\delta : \mathcal{C}^1(V,Z)\longrightarrow \mathcal{Z}^2(V,Z)$ such that
$$
\delta(g)(v_1,v_2)=g(v_1+v_2) + g(v_1) + g(v_2)\qquad \hbox{and}\qquad g(0_V)=0_Z.
$$ for all $v_1, v_2 \in V$.
Let $$\mathcal{B}^2(V,Z)=\delta(\mathcal{C}^1(V,Z))$$ and $$\mathcal{H}^2(V,Z)=\mathcal{Z}^2(V,Z)/\mathcal{B}^2(V,Z).$$

Central extensions $L_1$ and $L_2$ are called {\it equivalent} precisely if, there is an
isomorphism $\phi:L_{f_1}\longrightarrow L_{f_2}$ such that $\phi(v,\ast)=(v, \ast)$ if $v\in V$ and $\phi(\ast, z)=(\ast, z+\lambda(\ast))$ where $z, \lambda(\ast)\in Z$.

Any two equivalent extensions $L_1$ and $L_2$ are clearly isomorphic, but the converse is not true in general (for an example see the proof of Theorem \ref{N}).

\bl
Central extensions $L_{f_1}$ and $L_{f_2}$ corresponding to cocycles
$f_1$ and $f_2$ are equivalent if and only if $f_1 = f_2$ in $\mathcal{H}^2(V,Z)$.
\el

\bp
The map $\varphi=(\varphi_1,\varphi_2):L_{f_1} \longrightarrow L_{f_2}$,
with $\varphi_1(v,z)=v$ and $\varphi_2(v,z)=z+g(v)$, determines
an isomorphism if and only if $f_1(v_1,v_2)=f_2(v_1,v_2)+g(v_1 + v_2) + g(v_1) + g(v_2)$, i.e.,
$f_1=f_2$ in $\mathcal{H}^2(V,Z)$. This is because
$$
\varphi((v_1,z_1)\circ (v_2,z_2))=(v_1 + v_2, f_1(v_1,v_2) + z_1 + z_2 + g(v_1 + v_2))=
$$
$$
(v_1 + v_2, f_2(v_1,v_2) + z_1 + z_2 + g(v_1) + g(v_2))=\varphi(v_1,z_1)\circ \varphi(v_2,z_2).
$$
\ep

Let $\{v_1,...,v_n\}$ be a basis of $V_{\mathbf{F}_2}$; as before, we can identify $V$
with $P_n$ - the set of all subsets of $I_n$. The set $P_n$ has an ordering: $\sigma > \tau$ if
$|\sigma| > |\tau|$ or $|\sigma| = |\tau|$, $\sigma=\{i_1 < ... < i_k\}$, $\tau = \{j_1 < ... < j_k\}$ with $i_1=j_1, ... , i_s=j_s, i_{s+1} > j_{s+1}$.

Consider a subset $\mathcal{Z}^2_0(V,Z)\subset \mathcal{Z}^2(V,Z)$, where $f\in \mathcal{Z}^2_0(V,Z)$ if and only if
$f(\sigma, \{i\})=0$, $\{i\} \geq \max(\sigma)$, $\sigma\in P_n=V$. In what follows $\triangle$ stands for the set-theoretical difference.

\bl
$\mathcal{Z}^2(V,Z)=\mathcal{Z}^2_0(V,Z)\oplus \mathcal{B}^2(V,Z)$.
\el

\bp
First, consider the case when $f\in \mathcal{Z}^2_0(V,Z) \cap \mathcal{B}^2(V,Z)$. Then
$f=\delta(g)$, and for any $\sigma \in P_n$ and $\{i\}$ such that $\{i\} > \max(\sigma)$ we have
$$
f(\sigma, \{i\})= g(\sigma \cup \{i\}) + g(\sigma) + g(\{i\})=0.
$$
Then $g(\sigma)=\sum_{i \in \sigma} g(\{i\})$. Hence,
$$f(\sigma, \tau)= g(\sigma \triangle \tau) + g(\sigma) + g(\tau)=
\sum_{i \in \sigma \triangle \tau} g(\{i\}) + \sum_{i\in \sigma} g(\{i\}) + \sum_{i\in \tau}g(\{i\})
$$
$$
=\sum_{i\in \sigma \setminus \tau}g(\{i\}) + \sum_{i\in \tau \setminus \sigma}g(\{i\}) +
\sum_{i\in \sigma \cap \tau}g(\{i\})
$$
$$
+ \sum_{i\in \sigma \setminus \tau}g(\{i\}) +
\sum_{i\in \tau \cap \sigma}g(\{i\}) + \sum_{i\in \tau \setminus \sigma}g(\{i\})=0.
$$

Now, suppose $f\in \mathcal{Z}^2(V,Z)$. For $\sigma = (i_1,...,i_k)$ we define
$\sigma^s=(i_1,...,i_{s-1})$, $s>1$, and $g(\sigma)= \sum_{s=2}^{k}f(\sigma^s,\{i_s\})$,
assuming that $|\sigma|>1$ and $g(\{i\})=0$. Then $f+\delta(g)\in \mathcal{Z}^2_0(V,Z)$. Indeed, if
$\{i\}=\{i_{k+1}\}>\{i_k\}=\max(\sigma)$ then
$$(f+\delta(g))(\sigma,\{i\})=f(\sigma,\{i\})+
g(\sigma \cup \{i\}) + g(\sigma) + g(\{i\})$$
$$ = f(\sigma, \{i\}) + \sum_{s=2}^{k+1}f(\sigma^s,\{i_s\}) +
\sum_{s=2}^{k}f(\sigma^s,\{i_s\}) = 0,$$
as $\sigma^{k+1}=\sigma$ and $\{i_{k+1}\}=\{i\}$.
This yields that $f + \delta(g)\in \mathcal{Z}^2_0(V,Z)$
completing the proof of the lemma.
\ep

We call a pair $(\sigma, \tau)$ {\it regular}
if and only if $\sigma \triangle \tau > \sigma > \tau$. Note, that if $\emptyset \neq \sigma \neq \tau \neq \emptyset$ then
precisely one of the pairs $(\sigma,\sigma \triangle \tau)$, $(\sigma,\tau)$,
$(\sigma \triangle \tau,\tau)$, $(\sigma \triangle \tau, \sigma)$, $(\tau, \sigma)$, $(\tau, \sigma \triangle \tau)$ is regular.
A regular pair is called {\it strongly regular}
if $|\sigma| \geq |\tau| > 1$ or $|\sigma| \geq |\tau| = 1$ but $\{i\}< \max(\sigma)$,
where $\tau = \{i\}$ and $|\sigma|>1$.

\bl\label{NSRP}
The number of elements of the set of all strongly regular pairs is
$$\frac{1}{3}(2^{2n-1}+1) - 3\cdot 2^{n-1} + n + 1.$$
\el

\bp
Let $P$ be the set of all ordered pairs $(\sigma,\tau)$, $\sigma > \tau\neq \emptyset$.
Then $|P|= C_m^2, m=2^n-1$ and hence, the number of regular pairs is $\frac{1}{3}C_m^2=\frac{1}{3}(2^n - 1)(2^{n-1} - 1).$

If $(\sigma,\tau)$ is regular but not strongly regular then $\tau =\{i\}$,
$i \geq max(\sigma)$ or $\sigma = \{j\}$, $j>i$. Hence, for given $i$ we have $(2^{i-1} + n - 2i)$ regular but not strongly
regular pairs. Then the number of strongly regular pairs equals
$$
\frac{1}{3}(2^n - 1)(2^{n-1} - 1) - \sum_{i=2}^n (2^{i-1} + n -2i)
$$
$$ = \frac{1}{3}(2^n - 1)(2^{n-1} - 1) - 2^n -n^2 + n(n+1) + 1
$$
$$= \frac{1}{3}(2^{2n-1} + 1) - 3\cdot 2^{n-1} + n + 1.
$$
\ep

\bt\label{Basis}
The union of sets
$$\beal
\big\{(\{i\},\sigma\setminus \{i\},\tau)\; \big|\; (\sigma,\tau) \, \hbox{or} \, (\tau, \sigma) \,\, \hbox{strongly regular}, \\ \, \\ 
\qquad\qquad\qquad\sigma \cap \tau = \emptyset, \, \{i\} = max(\sigma \cup \tau) \in \sigma\big\}
\ena$$
and
$$
\big\{(\{j\},\mu,\lambda)\, \big| \, (\mu,\lambda) \,\, \hbox{strongly regular}, \,
\mu\cap\lambda\neq \emptyset, \, \{j\} = max(\mu \cap \lambda)\big\}.
$$
is a basis of the $\mathbf{F}_2$-space $S_1({\tt X})/S_2({\tt X})$,

Moreover,
$$
\dim_{\mathbf{F}_2}(S_1({\tt X})/S_2({\tt X})) =
\frac{1}{3}(2^{2n-1} + 1) - 3\cdot 2^{n-1} + n + 1,
$$
where $n := |{\tt X}|$.
\et

\bp
Let $F({\tt X})$ be a free 2-step nilpotent Steiner loop with free generators ${\tt X}=\{x_1, ..., x_n\}$.
$F({\tt X})$ can be realized as a central extension $L_f$ on $V\times Z$ for some Steiner loop cocycle
$f\in \mathcal{Z}^2_0(V,Z)$.

For any elements $(\sigma, s), (\mu, m), (\tau, t) \in L_f$ we have
$$\beal
\big((\sigma,s), (\mu,m), (\tau,t)\big) = \big[\big((\sigma,s)\circ(\mu,m)\big)\circ (\tau,t)\big] \circ
\big[(\sigma,s)\circ \big((\mu,m)\circ(\tau,t)\big)\big]\\ \, \\
\qquad =\big[(\sigma\triangle\mu, f(\sigma, \mu) + s + m)\circ(\tau,t)\big]\circ\big[(\sigma, s)\circ(\mu\triangle\tau, f(\mu, \tau) + m + t)\big] \\ \, \\
\qquad=(\sigma\triangle\mu\triangle \tau, f(\sigma, \mu)+f(\sigma\triangle\mu,\tau) + s+m+t)\\ \, \\
\qquad\qquad\circ(\sigma\triangle\mu\triangle\tau, f(\mu, \tau)+f(\sigma,\mu\triangle\tau) + s+m+t)\\ \, \\
\qquad= (\sigma\triangle\mu\triangle\tau\triangle \sigma\triangle\mu\triangle\tau, f(\sigma, \mu) + f(\sigma\triangle \mu, \tau) +
f(\mu, \tau) + f(\sigma, \mu \triangle\tau))\\ \, \\
\qquad = (\emptyset,  f(\sigma, \mu) + f(\sigma\triangle \mu, \tau) +
f(\mu, \tau) + f(\sigma, \mu \triangle\tau)).
\ena$$
Taking $s=m=t=0$ and identifying $(\sigma,0), (\mu,0)$ and $(\tau, 0)$ with $\sigma, \mu$ and $\tau$, respectively, we obtain the following relation involving associators:
\begin{equation}\label{formula}
(\sigma, \mu, \tau)=f(\sigma, \mu) + f(\mu, \tau) + f(\sigma \triangle \mu, \tau) +
f(\sigma, \mu \triangle \tau).
\end{equation}

Set $Z_f=f(V,V)$. Then by \eqref{formula} we get that $Ass(F({\tt X}))\subseteq Z_f$.

Next, we show that $Z_f\subseteq Ass(F({\tt X}))$. Let $\sigma, \tau \in V$ be such that $\sigma > \tau$.
If the pair $(\sigma, \tau)$ is not regular, then $\sigma > \sigma \triangle \tau$ and $f(\sigma, \tau)=f(\sigma \triangle \tau, \sigma)$
by the properties of Steiner loop cocycles. Note that the pair $(\sigma \triangle \tau, \sigma)$ is already regular.
Now, suppose that $(\sigma, \tau)$ is regular but not strongly regular. Then $\tau=\{i\}$ and $\{i\} \geq max (\sigma)$ or $\sigma=\{j\}$, $\tau=\{i\}$ and $j>i$. In this case
$f(\sigma, \{i\})=0$ by the definition of $\mathcal{Z}^2_0(V,Z)$. This means that it is enough to show that
$f(\sigma, \tau)\in Ass(F({\tt X}))$ for any strongly regular pair $(\sigma, \tau)$.

Let $(\sigma, \tau)$ be a strongly regular pair, that is, $|\sigma|\geq |\tau| > 1$ or $|\sigma|\geq |\tau|=1$ but
$\{i\}< max(\sigma)$, $\tau=\{i\}$ and $|\sigma|>1$. Furthermore, let $\{i\} = max(\sigma \cup \tau)$. In what follows
we use induction in $r:= |\sigma| + |\tau|$. Assume that $f(\sigma,\tau)\in Ass(F({\tt X}))$ if $(\sigma, \tau)$ is strongly regular and $|\sigma|+ |\tau|<r$. Consider the case where the pair $(\sigma, \tau)$ is strongly regular and $|\sigma|+ |\tau|=r$.

If $\sigma \cap \tau = \emptyset$ and $\{i\}\in \sigma$, then by \eqref{formula} we have:
$$\beal
(\{i\}, \sigma\setminus \{i\}, \tau)=f(\{i\}, \sigma\setminus \{i\}) + f(\sigma\setminus \{i\}, \tau)\\ \, \\
\qquad\qquad\qquad\quad + f(\{i\}\triangle (\sigma\setminus \{i\}),\tau) + f(\{i\},(\sigma\setminus \{i\}) \triangle \tau)
\ena$$
and
\begin{equation}\label{FSimTau1}
\beal
f(\sigma, \tau) = (\{i\}, \sigma\setminus \{i\}, \tau) + f(\{i\}, \sigma\setminus \{i\}) + f(\sigma\setminus \{i\}, \tau) \\ \, \\
\qquad\qquad+ f(\{i\}, (\sigma\setminus \{i\})\triangle \tau)=(\{i\}, \sigma\setminus \{i\}, \tau) + f(\sigma\setminus \{i\}, \tau),
\ena
\end{equation}
since $f(\{i\},\sigma\setminus \{i\})=0$ and $f(\{i\}, (\sigma\setminus \{i\})\triangle \tau)=0$. By the induction assumption, $f(\sigma\setminus \{i\}, \tau)\in Ass(F({\tt X}))$
as $|\sigma\setminus \{i\}|+|\tau|<|\sigma|+|\tau|$ and hence $f(\sigma, \tau)\in Ass(F({\tt X}))$. Similarly, we can prove the same fact in the case where
$\sigma \cap \tau =\emptyset$ and $\{i\}\in \tau$.

If $\sigma \cap \tau \neq \emptyset$ and $\{j\}=max(\sigma \cap \tau)$, then by \eqref{formula} we obtain that
$$
(\{j\}, \sigma, \tau)=f(\{j\},\sigma) + f(\sigma,\tau) + f(\{j\}\triangle \sigma, \tau) + f(\{j\}, \sigma \triangle \tau)
$$
and
\begin{equation}\label{FSigmaTau2}
f(\sigma, \tau) = (\{j\}, \sigma, \tau) + f(\{j\}, \sigma) + f(\sigma\setminus \{j\}, \tau) + f(\{j\}, \sigma \triangle \tau).
\end{equation}
Each summand in the right-hand side
is contained in $Ass(F({\tt X}))$. Namely, as $f$ is a Steiner loop cocycle, we have $f(\{j\},\sigma)=f(\sigma,\{j\})=f(\sigma\setminus \{j\}, \{j\})$. Then $|\sigma\setminus \{j\}|+|\{j\}|=|\sigma|<|\sigma|+|\tau|$ yields by the induction hypothesis that $f(\{j\},\sigma)\in Ass(F({\tt X}))$. Similarly, we have that
$f(\sigma\setminus \{j\}, \tau)\in Ass(F({\tt X}))$ by induction, since $\{j\}\in \sigma$ and $|\sigma\setminus \{j\}|+|\tau|<|\sigma|+|\tau|$. Also, $f(\{j\}, \sigma\triangle \tau)\in Ass(F({\tt X}))$ by the induction, since $|\{j\}|+|\sigma\triangle \tau|=1+|\sigma| + |\tau| -2|\sigma\cap \tau|<|\sigma| + |\tau|$ as $\sigma \cap \tau \neq \emptyset$. This implies that $f(\sigma, \tau)\in Ass(F({\tt X}))$.

Summarizing the above discussions, we get that $Z_f\subseteq Ass(F({\tt X}))$. Hence $Z_f = f(V,V) = Ass(F({\tt X}))=S_1({\tt X})/S_2({\tt X})$.

Moreover, by \eqref{FSimTau1} and \eqref{FSigmaTau2} we obtain that
$$\beal
f(V,V) \subseteq Span_{\mathbf{F}_2}\Big\{ (\{i\}, \sigma\setminus \{i\}, \tau), (\{j\}, \mu, \lambda) \\ \, \\
\qquad\qquad\quad \Big| (\sigma, \tau), (\mu, \lambda) \; \hbox{strongly regular},
\sigma\cap\tau=\emptyset, \{i\}=max(\sigma\cup\tau)\\ \, \\
\qquad\qquad\quad \hbox{and} \, \mu\cap\lambda\neq\emptyset, \{j\}=max(\mu\cap\lambda) \Big\}=:W.
\ena$$

Finally, we have that $f(V,V)\subseteq W\subseteq Ass(F({\tt X}))=f(V,V)$ which implies $Ass(F({\tt X}))=W$ and
$$
\beal
dim(S_1({\tt X})/S_2({\tt X}))=dim(f(V,V))=dim(Ass(F({\tt X})))=dim(W)\\ \, \\
\qquad\qquad=\big|\{\hbox{strongly regular pairs}\}\big|=\frac{1}{3}(2^{2n-1}+1) - 3\cdot 2^{n-1} + n + 1
\ena$$
where the last equality holds by Lemma \ref{NSRP}.
\ep

Note that there are exactly $80$ non-isomorphic Steiner triple systems of order $15$. Moreover, there is only one nilpotent non-associative Steiner loop $S_{16}$ of order $16$ (cf. \cite{GAP}), and it corresponds to the system N.2 in \cite{MPhR} p. 19. Furthermore, $S_{16}$ has the GAP id ${\tt SteinerLoop( 16, 2)}$; the label $2$ indicates the system order as in the list established in monograph \cite{CR}.
The Steiner loop $S_{16}$ is $3$-generated and has the nilpotency class $2$. In what follows we describe all $3$-generated Steiner loops of nilpotency class $2$.

\bt\label{N}
There exist three non-isomorphic non-associative $3$-generated Steiner loops of nilpotency class $2$ and their orders are $16,$ $32$ and $64$, respectively.
\et

\bp
Let $S({\tt X}=\{x_1, x_2, x_3\})$ be the $3$-generated free Steiner loop of nilpotency class 2 and $Z=<(z_1, z_2, z_3)>$ be a center of $S({\tt X}).$  By Theorem \ref{Basis}, we can choose
$z_1=(x_1,x_2,x_3),$ $z_2=(x_2,x_1,x_3),$ $z_3=(x_3,x_2,x_1x_3).$

Let $G=AutS({\tt X})$ be a group of automorphisms of the loop $S({\tt X}).$ Since $S({\tt X})$
is free and $Z$ is a $G-$invariant subloop of $S({\tt X})$, we have an epimorphism $\phi:G\to GL_3({\bf F}_2)$ with $ker(\phi)=\{\rho \in G | x^{\rho }=xz, z\in Z\}.$

The group $G$ acts on the ${\bf F}_2-$space $Z$, and this action depends only on the images of the elements of $G$ in $GL_3({\bf F}_2).$ As the group $GL_3({\bf F}_2)$ is simple and has no
non-trivial $2-$dimensional ${\bf F}_2-$representations, the $G-$module $Z$ is irreducible.

If $P$ is a $3$-generated Steiner loop of nilpotency class $2$, then there is a canonical
epimorphism $\pi:S({\tt X})\to P $ and $ker(\pi)\subseteq Z.$ Note that for any other $3$-generated Steiner loop $Q$ and canonical epimorphism $\psi:S({\tt X})\to Q$ the
loops $P$ and $Q$ are isomorphic if and only if there exists $\varrho\in GL_3({\bf F}_2)$
such that $ker(\pi)^{\varrho}=ker(\psi).$ Indeed, if $ker(\pi)^{\varrho}=ker(\psi)$ then
$\varrho$ induces an isomorphism $\bar{\varrho}:P=S({\tt X})/ker(\pi)\longrightarrow S({\tt X})/ker(\psi)=Q$. Conversely, let $\bar{\varrho}:P=S({\tt X})/ker(\pi)\longrightarrow S({\tt X})/ker(\psi)=Q$ be an isomorphism between $P$ and $Q$. Then $\bar{\varrho}$ induces a homomorphism
$\upsilon:S({\tt X})\longrightarrow S({\tt X})/ker(\psi)$ with $\upsilon=\bar{\varrho}\circ \pi$.
For every $x\in \tt X$ one can choose $\sigma(x)\in \tt X$ such that $\upsilon(x)=\sigma(x)ker(\psi)\in
S({\tt X})/ker(\psi)$. Since the loop $S({\tt X})$ is free, there is a unique homomorphism
$\bar{\upsilon}:S({\tt X})\longrightarrow S({\tt X})$ satisfying $\bar{\upsilon}(x)=\upsilon(x)$
for all $x\in \tt X$. It is easy to see that $\bar{\upsilon}\in Aut S({\tt X})$ and $\bar{\upsilon}(ker (\pi))=ker(\psi)$.

The group $GL_3({\bf F}_2)$ acts transitively on $Z\setminus \{1\}$ and on the set of the two dimensional ${\bf F}_2-$subspaces of $Z$. As $Z$ is a three dimensional irreducible
$GL_3({\bf F}_2)-$module, there exists a unique $3$-generated Steiner loop of nilpotency class $2$ for each of the orders $16$, $32$ and $64$.
\ep

\noindent
{\bf Acknowledgement}
\vskip 5pt
The first autor is grateful to FAPESP and CNPq, Brazil; to RFBR, Russia, grant 16-01-0577a (pp. 1-5) and to
Russian Science Foundation, project 16-11-10002 (pp. 6-10) for financial support. The second author has been supported by FAPESP Grant - 2015/17611-8. The third author has been supported by CNPq, process 307824/2016-0. The fourth author has been supported by FAPESP Grant - 11/51845-5, and expresses her gratitude to IMS, University of S\~{a}o Paulo and the Math. Dept., Penn State University, for the warm hospitality.

The authors are grateful to the referee for careful reading of the manuscript and a numerous
queries aiming at a substantial improvement of the paper.

\vskip 5pt

Alexander Grishkov\\
IMS, University of S\~{a}o Paulo 05508-090 S\~{a}o Paulo, SP, Brazil \vspace{1mm}\\
Omsk State University n.a. F.M. Dostoevskii, Pr. Mira 55-A, 644053 Omsk, Russia\\
{\it E-mail}: {\it {}grishkov@ime.usp.br}\\

Diana Rasskazova\\
IMS, University of S\~{a}o Paulo, 05508-090 S\~{a}o Paulo, SP, Brazil \vspace{1mm}\\
{\it E-mail}: {\it {}ivakirjan@gmail.com}\\

Marina Rasskazova\\
Omsk State Technical University, 644050 Omsk, Russia \vspace{1mm}\\
{\it E-mail}: {\it {}marinarasskazova1@gmail.com}\\

Izabella Stuhl\\
Math Dept, Penn State University, State College, PA 16802, USA\\
{\it E-mail}: {\it {}ius68@psu.edu}\\

\end{document}